\documentclass[12pt]{amsart}

\usepackage{amsmath}
\usepackage{amssymb}
\usepackage{latexsym}
\newenvironment{proof*}{\vskip 2mm\noindent {}}{$\hfill \Box$ \vskip 2mm}

\newtheorem{theorem}{Theorem} 
\newtheorem*{theorem*}{Theorem} 
\newtheorem{lemma}{Lemma} 
\newtheorem{corollary}{Corollary} 
\newtheorem{proposition}{Proposition} 
\newtheorem{definition}{Definition} 
 
\newtheorem{claim}{Claim }

                     \newcommand{\Dd}{{\mathbb{D}}} 
                      
                     \newcommand{\Zz}{{\mathbb{Z}}}

                     \newcommand{\calf}{{\mathcal{F}}}

                     \newcommand{\call}{{\mathcal{L}}}

                     \newcommand{\calp}{{\mathcal{P}}}
                     \newcommand{\calq}{{\mathcal{Q}}}
                     \newcommand{\cals}{{\mathcal{S}}}

\newcommand{\di}{{\mathbb D}}
\renewcommand{\a}{\alpha}
\newcommand{\al}{\alpha} 
\newcommand{\la}{\lambda}

\newcommand{\eit}{e^{i\theta}}
\newcommand{\emit}{e^{-i\theta}}

\newcommand{\Stolat}{\Gamma_\a(\eit)}

\newcommand{\disak}{(1-|a_k|)}

\newcommand{\sak}{\{a_k\}}

\newcommand{\svn}{\{v_n\}}
\newcommand{\sbk}{\{b_k\}}
\newcommand{\swn}{\{w_n\}}
\newcommand{\sumk}{\sum_{k=1}^\infty}

\renewcommand{\Re}{\mbox{Re}}

\newcommand{\tg}{{\tilde g}}

\begin{document} 

\pagestyle{plain}

\title{\bf Decrease of bounded holomorphic 
functions along discrete sets} 

                     \author{Jordi Pau 
 \ \ \ Pascal J. Thomas} 

\begin{abstract}
We provide results of
uniqueness for holomorphic functions in the Nevanlinna class bridging
those previously obtained by Hayman and Lyubarskii-Seip. Namely,
we propose certain classes of hyperbolically separated sequences
in the disk, in terms of the rate of non-tangential accumulation to the
boundary (the outer limits of this spectrum of classes being
respectively the sequences with
a non-tangential cluster set of positive measure,
and the sequences satisfying the Blaschke condition); and for each
of those classes, we give a critical condition of radial decrease
on the modulus which will force a Nevanlinna class function to vanish
identically.
 \end{abstract}

\maketitle


\footnotetext{The first  author was partially supported by DGES grant
PB98-0872 and CIRIT grant 2001 SGR00431, and both authors' travels were supported
by a grant in the framework of the Comunitat de Treball dels Pirineus.}


\vskip1cm 

\section{Definitions and previous results}
\label{defs}

Let $\mathbb D$ be the unit disk in the complex plane. We are interested in the
allowable decrease of the modulus of a 
nontrivial bounded holomorphic function 
$f$ along a discrete sequence $\sak$ in $\di$. 
Notice that those problems are actually
the same if we replace the class of bounded holomorphic functions,
$H^\infty$, by the Nevanlinna class $N$, since
on the one hand $H^\infty \subset N$ and on the other hand, for any
$f$ in $N$, we can write $f = f_1/f_2$ where $|f_2(z)| \le 1$ whenever $z \in \Dd$,
thus $|f_1(z)| \le |f(z)|$, and $f_1 \in H^\infty$. This also means 
in particular that 
any Hardy space $H^p$ could be substituted for $H^\infty$. On the other hand
it seems clear that the situation in Bergman spaces has to be quite different,
see \cite{Bo} for related results. Previous results about 
those problems can be found in \cite{Kha},
\cite{Ei1}, \cite{Ha}, \cite{Lu-Se} and recent work of
Eiderman and Ess\'en reported on in \cite{Ei2}. Results in a similar
spirit about positive harmonic
functions could already be found in \cite{Be}.

We now proceed to describe the results in \cite{Ha} and \cite{Lu-Se} that 
served as an inspiration for the present work.

Let $\Stolat$ stand for the Stolz angle of aperture $\alpha >0$ with vertex 
at the point $\eit$ on the circle,
$$
\Stolat := \{ \zeta \in \di : |1-\zeta \emit| < (1 + \a ) (1-|\zeta|) \} .
$$
To simplify notations, we sometimes use the single lower-case letters $a$, 
$b$ to stand for the sequences $\sak$, $\sbk$...
We denote by $NT(a)$ the set 
of points on the circle which are non-tangential limit points of the 
sequence $\sak$, 
\begin{equation}
\label{setNT}
NT(a) := \{ \eit \in \partial\di : \exists \a >0 : \eit \in 
\overline{\Stolat \cap \sak} \}.
\end{equation}

A first result about decrease of bounded holomorphic functions is given, 
along with more elaborate theorems, in \cite{Ha}:

\begin{theorem*}[A]
Any function $f$ in $H^\infty$
such that $\lim_{k\to\infty} f(a_k) =0$ must vanish identically 
if and only if $|NT(a)|>0$.
\end{theorem*}

When 
$\sak$ is a Blaschke sequence, that is, when it satisfies the Blaschke condition 
$$
\sumk \disak < \infty ,
$$
there exists a nontrivial bounded holomorphic function vanishing exactly 
on the sequence. It is easy to prove directly that a Blaschke sequence
$a$ verifies $|NT(a)|=0$ (see below). 
Thus questions of decrease only make sense for non-Blaschke sequences, and 
we are especially interested in the intermediate case of non-Blaschke sequence with $|NT(a)|=0$. 

It is often useful to restrict attention to 
sequences which are {\it separated\/}
in the Gleason (or pseudo-hyperbolic)
distance $d_G (z,w):=\left| \frac{z-w}{1-\bar w z}
\right|$,  that is, 
$$
\delta := \inf_{j \neq k} 
\left| \frac{a_k - a_j}{1-\bar a_k a_j} \right| > 0 .
$$
The number $\delta$ is sometimes called {\it separation constant\/}.

Following the idea of \cite{Lu-Se}, given a class of sequences $\cals$ in the disk, we shall say:

\begin{definition}
\label{essmin}
A nonincreasing function $g$ from $[0,1)$ to $(0,1]$,
tending to $0$ as $x$ tends to $1$,
is an ($H^\infty$-){\em essential minorant\/} for $\cals$ 
if and only if given any
 sequence $a=\sak$ in $\cals$, any bounded holomorphic function $f$
verifying $|f(a_k)| \le g(|a_k|)$ for all $k$
must vanish identically. 
\end{definition}

For instance, a consequence of the ``if" part of Theorem A is that any 
function $g$ such that $\lim_{x\to1} g(x)=0$ 
is an essential minorant for the class of sequences $a$
such that $|NT(a)|>0$. 

This is the main result of \cite{Lu-Se}:

\begin{theorem*}[B]
A nonincreasing function $g$  
is an $H^\infty$-essential minorant for the separated 
non-Blaschke sequences if and only if
\begin{equation}
\label{critdecr}
\int_0^1 \frac{dr}{(1-r)\log \left(\frac1{g(r)}\right)} < \infty.
\end{equation}
\end{theorem*}

This theorem reveals the following remarkable fact: 
given a function $f$ in $H^\infty$, when one measures the decrease 
of $|f|$ 
outside of a fixed hyperbolic neighborhood of the zero set of $f$, it 
cannot go below a certain critical velocity given by condition 
(\ref{critdecr}).

We introduce a transformation of the function $g$ measuring the decrease,
which is related to the fact that we will work with the (sub-)harmonic 
function $\log |f|$, and that the sequences are separated.

\begin{definition}
\label{fcng}
For any nonincreasing function
$g$ from $[0,1)$ to $(0,1]$, tending to $0$ as $x$ tends to $1$,
we write
$$
\tilde g (\la) := \log \frac1{g(1-2^{-\la})} , \quad \la \ge 0.
$$
\end{definition}

It is then elementary to see that condition (\ref{critdecr}) in Theorem B can 
 be restated as :
\begin{equation}
\label{discretedecr}
\sum_0^\infty \frac{1}{\tilde g(n)} < \infty .
\end{equation}

\section{The main theorem}
\label{mainthm}

We aim to bridge the gap between Theorems A and B by introducing classes 
of sequences which mediate between those involved in each of those results.

Notice first that one could define a set $NT_\a(a)$ by fixing the number $\alpha$ 
in (\ref{setNT}). Two such choices would differ from each other (and thus 
from the whole $NT(a)$) only by a set of measure $0$ \cite{Th}. Hence we 
will neglect sets of linear measure $0$, fix a value of $\a$, and write 
$NT(a)$ for what is really $NT_\a(a)$.

In order to measure the density of a sequence $a=\sak$ as it accumulates to 
the unit circle, we will consider the following function on the circle.
$$
\phi_a (\eit) := \# \sak \cap \Stolat .
$$
This is the same function which was denoted by $\Gamma_\gamma$ in 
\cite[eqn. 0.1]{NPT}. One sees easily that $\sak$ satisfies the Blaschke 
condition if and only if $\phi_a \in L^1 (\partial \di)$; and that the set 
$NT(a)$ is the set of $\eit$ such 
that $\phi_a (\eit) = \infty$ \cite{NPT}, \cite{Th}. Thus if $a$ is 
Blaschke, $|NT(a)|=0$.

Because of the rotation invariance of the essential minorants we are considering,
most of the relevant information about $\phi_a$ is given by
$$
m_a (n) := \left| \{ \zeta \in \partial \di : \phi_a (\zeta) \ge n
\} \right| , \ n \in \Zz_+ .
$$
For instance, $\int_{\partial \Dd} \phi_a = \sum_{n \ge 1} m_a (n)$ and
$|\{ \phi_a = \infty \}|= \lim_{n\to\infty} m_a (n)
= \inf_n  m_a (n) = |NT(a)|$.

We will now introduce classes of non-Blaschke sequences, the size of which 
is measured by the behaviour of their associated function $\phi_a$.

\begin{definition}
\label{classes}
Let $\svn$, $\swn$ be nonincreasing
bounded sequences with nonnegative values such that 
$\sum_0^\infty v_n = \infty$,
$\sum_0^\infty w_n = \infty$. Define
$$
\cals_w := \left\{ a \mbox{ separated} : \sum_0^\infty m_a (n) w_n = \infty
\right\} ,
$$
and
$$
\call_v := \left\{ a \mbox{ separated} : \liminf_{n \to \infty} m_a (n) / v_n > 0
\right\} .
$$
\end{definition}
Clearly, if there is a $C >0$ such that $w \le C w'$, then
$\cals_{w} \subset \cals_{w'}$ and $\call_{w'} \subset \call_{w}$.
Note that the class of all separated non-Blaschke sequences coincides with $\cals_1$,
where $1$ stands for the constant sequence $w_n=1$ for all $n$, 
and the set of separated sequences with 
$|NT(a)|>0$ coincides with $\call_1$. 

Also, for any $v$, $w$ satisfying the
conditions above, $\call_1 \subset \call_v \subset \cals_1$, and
$\call_1 \subset \cals_w \subset \cals_1$, in particular, the classes 
we have defined are never empty and never contain any Blaschke sequence. 
Finally, it is also easy to see that $\call_1 = \bigcap \call_v = \bigcap 
\cals_w$, and $\cals_1 = \bigcup \call_v = \bigcup \cals_w$. 

More generally, we have the following properties:
\begin{enumerate}
\item
$\call_v \subset \cals_w$ if and only if $\sum_0^\infty v_n w_n =
\infty$. 
\item
$\cals_w = \bigcup \{ \call_v : \sum_0^\infty v_n w_n =
\infty \}$.
\end{enumerate}

Our main result is 
\begin{theorem}
\label{main}
\begin{enumerate}
\item 
\label{sum}
$g$ is an essential minorant for the class $\cals_w$ if and only 
if 
$$
\sum_0^\infty \frac{w_n}{\tilde g(n)} < \infty ;
$$
\item 
\label{liminf}
$g$ is an essential minorant for the class $\call_v$ if and only 
if 
for any subset $E$ of $\Zz_+$ such that 
$\sum_{n\in E}v_n < \infty$, and
for any positive integer $C$,
$$\limsup_{n \notin E} \tilde g([\frac{n}C]) v_{n} 
= \infty .
$$
\end{enumerate}
\end{theorem}
Here $[x]$ denotes as usual the largest integer which is smaller
or equal than the real number $x$.

Part (1) with $w_n = 1$ gives Theorem B \cite{Lu-Se}.
On the other hand, Part (2) shows that given any $g$ satisfying
the hypotheses in Definition \ref{essmin}, we can find a class of
sequences strictly greater than $\call_1$ such that $g$ is still
an essential minorant. This sharpens Theorem A (for separated sequences only).

To take into account zero-sets of bounded analytic functions,
we will have
 to introduce the following
modified versions of the classes $\call_v$.

\begin{definition}
\label{Lmodif}
\[
\call'_v := \left\{ a : \exists b\mbox{, Blaschke sequence,
such that } a \cup b \in \call_v
\right\} .
\]
\end{definition}

Note that this makes the class stable under removal of a Blaschke sequence:
if $a \in \call'_v$ and $b$ is a Blaschke sequence, then $a\setminus b 
\in \call'_v$.
In many cases, this makes no difference.

\begin{lemma}
\label{Lstable}
If $\svn$ verifies that there exists $\eta_1 > 1, \eta_2 > 0$
such that $v_{[\eta_1 n]} \ge \eta_2 v_n$ (which implies
$v_n \ge \frac{C}n$ for all $n \ge 1$), then 
$\call_v =  \call'_v$.
\end{lemma}

When this is the case, the unpleasantly complicated
statement of Theorem \ref{main} (\ref{liminf}) can be 
simplified.

\begin{corollary}
\label{cor}
Under the hypotheses of Lemma \ref{Lstable}, 
$g$ is an essential minorant for the class $\call_v$ if and only 
if 
$$
\limsup_{n \to \infty} \tilde g(n) v_{n} 
= \infty .
$$
\end{corollary}

\section{Reduction to the zero-free case}

It will be useful to have some preliminary result on the stability of the 
classes of sequences we have defined above.

\begin{lemma}
\label{union}
Suppose that $a$, $b \notin \cals_w$. Then $a \cup b \notin \cals_w$.

In particular, suppose we remove a Blaschke sequence $b$ from $a \in 
\cals_w$. Then $a \setminus b \in \cals_w$.
\end{lemma}

\begin{proof}
We have $\phi_{a\cup b} = \phi_{a}+ \phi_{b} $, so
$$
\{ \zeta \in \partial \di : \phi_{a\cup b} (\zeta) \ge \la \}
\subset
\{ \zeta \in \partial \di : \phi_a (\zeta) \ge \frac\la2
\}
\cup
\{ \zeta \in \partial \di : \phi_a (\zeta) \ge \frac\la2
\}
$$
therefore $m_{a\cup b} (n) \le m_a ([n/2]) + m_b ([n/2])$, and by the 
monotonocity of $w$,
\begin{eqnarray*}
\sum_0^\infty m_{a\cup b} (n) w_n
&\le&
\sum_0^\infty m_{a} ([n/2]) w_{[n/2]}
+
\sum_0^\infty m_{b} ([n/2]) w_{[n/2]} \\
&=& 2 \sum_0^\infty m_{a} (n) w_n + 2  \sum_0^\infty m_{b} (n) w_n
< \infty.
\end{eqnarray*}
\end{proof}

The next step is to reduce our problem to one  about
positive harmonic functions. There is no loss of generality in assuming
always that the function $g$ is bounded above by $1$.

\begin{lemma}
\label{harm}
A function $g$ is an essential minorant for the class $\cals_w$
(resp. $\call_v$)
if and only if for any  sequence $a$ in $\cals_w$ (resp. $\call'_v$)
there exist no (positive) harmonic function $h$ on the disk such that 
$$
h (a_k) \ge \log \frac1{g(|a_k|)}, \mbox{ whenever }k \in \Zz_+.
$$
\end{lemma}

\begin{proof}
One direction is clear : if for some $a$ in $\cals_w$,
a harmonic function $h$ as above exists, then the 
function $f:= \exp (-h -i \tilde h)$, where $\tilde h$ denotes
the Hilbert transform of $h$, is holomorphic and bounded, 
and $|f(a_k)| \le g(|a_k|)$, 
so that $g$ is not an essential minorant for $\cals_w$.

If on the other hand $a \in\call'_v$, choose a Blaschke sequence
$b$ so that  $a \cup b \in\call_v$, denote by $B$ the Blaschke product
with zeros on the sequence $b$, then
 $f:= B \exp (-h -i \tilde h)$ will verify $|f(a_k)| \le g(|a_k|)$.

The proof of the converse
 is essentially the same as that of Lemma 4 in \cite[p. 123]{NPT}.

Suppose that $g$ is not an essential minorant. Then there
exist $a'$ in $\cals_w$ (resp. $\call_v$) and
$f$ in $H^\infty \setminus \{0\}$ 
such that $\|f\|_\infty \le 1$ and
 $|f(a'_k)| \le g(|a'_k|)$ for all $k$.
Then $f=B f_1$, where $f_1$ is zero-free and $B$ is a Blaschke product,
\[
B(z) = z^m \prod_k \frac{|b_k|}{b_k} \frac{b_k-z}{1-\bar b_k z},
\]
where $b_k \in \di$, $\sum_k (1-|b_k|) < \infty$. Let
$\delta=\inf\{d(a'_k, a'_j)\colon k\not = j\}> 0$,
by the separation condition.
First notice that 
for each $k$, there is at most
one point
$a'_j$ such that $d_G (b_k, a'_j) < \delta/2$, so the
sequence
$$
b' := \left\{ a'_j : d_G (\{b_k\}, a'_j) < \delta/2 \right\}
$$
is a Blaschke sequence. Define a sequence $\sak$ by \ $a := a' \setminus
b'$. 
If $a' \in \call_v$, by construction
 $a \in \call'_v$. If $a' \in \cals_w$, by Lemma \ref{union}, 
$a \in \cals_w$. 

In \cite{NPT} on page 124, 
lines 3 to 17, it is proved that there exists 
 a holomorphic function $f_2$ in the unit disk such that
$|e^{f_2}(z) | \le |B(z)|$ for any $z$ such that
$d (z, B^{-1}(0)) \ge \delta/2$.

Then the function $h := -\Re \log f_1 - \Re f_2$
satisfies the conclusion of Lemma \ref{harm}.
\end{proof}

\section{A dyadic partition}
\label{dyadgene}

Consider the following partition of $\partial \di$ in dyadic arcs,
for any $n$ in $\Zz_+$:
\[
I_{n,k} := \{ \eit : \theta \in [{2\pi}k 2^{-n},{2\pi}(k+1) 2^{-n})
\}, \ 0 \le k < 2^n .
\]
To this we associate the "dyadic cubes" 
\[
Q_{n,k} := \{ r \eit : \eit \in I_{n,k} , 1- 2^{-n} \le r < 1- 2^{-n-1} \}.
\]
The following property is well known and easy to check.

\begin{equation}
\label{dyadsize}
\exists \delta_0 < 1 \mbox{ s.t. } \forall z_1,z_2 \in
Q_{n,k} , \quad d_G(z_1,z_2) \le \delta_0 .
\end{equation}

Whenever we consider a sequence, it will be natural to consider
the set $\calq$ of the dyadic cubes which meet the sequence. On the other hand,
from equation (\ref{dyadsize}) and
the fact that surface area is doubling with respect to the 
Gleason distance, we deduce
that if a sequence $\sak$ is separated, then there exists
an integer $N = N(\delta)$, where $\delta$ is the separation constant, 
so that 
\begin{equation}
\label{pointsinbox}
\# \sak \cap Q_{n,j} \leq N\mbox{, whenever }n \ge 0, 0\le j < 2^n .
\end{equation}

A Stolz angle could be defined as the union of all hyperbolic balls of a fixed 
size that are centered on a given radius (or that intersect it). Analogous 
facts hold for dyadic cubes. 
For any given $M$ in $\Zz_+$, let 
\[
Q_{n,j}^M := \bigcup_{l:j-M \le l \le j+M} Q_{n,l} 
\]
(where the integers $j$ and $l$ are taken mod $2^n$).
The proof of the following fact is best left to the reader.

\begin{lemma}
\label{neighborscover}
For any positive number $\alpha$, there exists $M_1=M_1(\al)$ such that 
\[ 
\Gamma_\al (\eit) \subset \bigcup \left\{Q_{n,j}^{M_1} :  Q_{n,j} \cap 
\{r\eit, 0\le r <1\} \neq \emptyset \right\} .
\]
\end{lemma}

\section{Proof of the sufficient condition}
\label{proofsuf}

We begin with a simple remark. Let $p$ be any nonnegative function on the disk. 
The 
 nontangential maximal function associated to $p$
is 
$$
p^* (\eit) := \sup_{\Stolat} p .
$$

From the fact that if $h$ a positive
harmonic function, then $h^*$ is weak $L^1(\partial \Dd)$
(see e.g. \cite[Theorem 5.1, p. 28]{Ga}), we deduce :

\begin{lemma}
\label{maxfcn}
If $h$ is harmonic in the disk and $h(z) \ge p(z) \ge 0$,
whenever $z\in \Dd$, then
$p^*$ is weak $L^1 (\partial \Dd)$, that is,
for any positive number $\lambda $, $|\{ p^* > \lambda \}| \le C/ \lambda$.
\end{lemma}

\begin{proof*}{\it Proof of sufficiency in Theorem \ref{main} (\ref{sum}) :}

Suppose that we have $f \in H^\infty \setminus \{0\}$,
a function $g$ as in Definition \ref{fcng}
such that $\sum_n w_n/\tg(n) < \infty$ and $a' \in \cals_w$ such that
$|f(a'_k) | \le g(|a'_k|)$. Then by Lemma \ref{harm}, we may assume that
we have in fact $h(a_k)  \ge \log (1/ g(|a_k|))$, with $h$ a harmonic
function and $a$ in $\cals_w$. 

We set $p(z):= \log \frac1{g(|a_k|)} \ge 0$ when $z= a_k$, $p(z) = 0$ otherwise.
We now want to relate the sets where 
$p^*$ is large  to those where $\phi_a$ is large.
 
\begin{lemma}
\label{phip}
There exist a constant $M \ge 1$ (depending only on $\alpha$ and the
constant $N$ in the hypotheses above) such that
\begin{equation}
\label{setincl}
\{ \phi_a \ge n \} \subset \{ p^* \ge \tg([\frac{n}M]) \}.
\end{equation}
\end{lemma}

\begin{proof}
We prove the contraposite inclusion.
Suppose that $p^* (\eit) < \tg(m)$. Then, for any $a_k$ in $\Stolat$, 
$\log \frac1{g(|a_k|)} < \tg(m)$, that is to say, $|a_k| < 1-2^{-m}$.
Then, by Lemma \ref{neighborscover}, we see that $a \cap \Stolat$ is contained
in the union of sets $Q_{n,j_n}^{M_1}$ for $0 \le n < m$, and therefore, 
because of (\ref{pointsinbox}),
$$
\phi_a (\eit) =
\# (a \cap \Stolat) \le N  (2M_1 (\al) + 1) m =: M(N, \al) m.
$$
Choosing $m$ approximately equal to $n/M(N, \al)$, we obtain 
that $\phi_a (\eit)<n$, therefore, with a slightly 
modified $M$, $\{ p^* < \tg([\frac{n}M]) \} \subset \{ \phi_a < n \}$.
\end{proof}

By Lemma \ref{maxfcn}, which applies because of the assumptions made
at the start of the proof of Theorem \ref{main} (1), the inclusion
(\ref{setincl}) implies that 
\begin{equation}
\label{domination}
m_a(n) \le \frac{C}{\tg([\frac{n}M])}.
\end{equation}
 Then, by the monotonicity of
$\{w_n\}$,
\[
\sum_n w_n m_a(n) \le C \sum_n \frac{w_n}{ \tg([\frac{n}M])}
\le  C \sum_n \frac{w_{[\frac{n}M]}}{ \tg([\frac{n}M])}
\le  CM \sum_n \frac{w_n}{ \tg(n)} < \infty,
\]
which contradicts the fact that $a \in \cals_w$.
\end{proof*}

\begin{proof*}{\it Proof of sufficiency in Theorem \ref{main} 
(\ref{liminf}) :}

Suppose that we have $f \in H^\infty \setminus \{0\}$,
a function $g$ as in Definition \ref{fcng}
and $a' \in \call_v$ such that
$|f(a'_k) | \le g(|a'_k|)$. Then by Lemma \ref{harm}, we may assume that
we have in fact $h(a_k)  \ge \log (1/ g(|a_k|))$, with $h$ a harmonic
function and $a \in \call'_v$. 
 Lemma \ref{phip} 
and its consequence (\ref{domination}) then apply to the
separated sequence $a$ and the decreasing function $g$.

To get a contradiction, we assume that
for any positive integer $M$ and any subset $E$ of $\Zz_+$ such that
$\sum_{n\in E} v_n <\infty$,
$\limsup_{n\notin E} \tg([\frac{n}M]) v_{n} = \infty$. 

Let $b$ be a Blaschke sequence such that $m_{a\cup b}(n) \ge c_0 v_n$
for $n$ large enough.
As in the proof of Lemma \ref{union}, 
$m_{a\cup b} (n) \le m_a ([n/2]) + m_b ([n/2])$. Let $E$ be the set of
indices $n$ such that $m_b ([n/2]) \ge \frac12 m_{a\cup b} (n)$. Since
by hypothesis $\sum_n m_b(n) < \infty$, we do have $\sum_{n \in E} v_n 
< \infty$.

For $n \notin E$, 
\[
m_a ([n/2]) \ge \frac12 m_{a\cup b} (n)  \ge \frac{c_0}2 v_n.
\]
Therefore, applying (\ref{domination}),
\[
\tg([\frac{n}M]) \le \frac{C}{m_a ([n/2])} \le \frac{C'}{v_n},
\]
which contradicts our assumption.
\end{proof*}

\section{Proof of the necessary condition}
\label{proofnec}

\begin{lemma}
\label{mucantor}Given a function $g$ as in Definition \ref{fcng} such that $\sum_ n
1/\tilde{g}(n)=\infty$, there exist a constant
$C>0$, a function $h$ harmonic in the disk, a non-empty sequence $p$ and a
Blaschke sequence $b$ such that 
\begin{enumerate}
\item
For any $p_k$ in $p$, $h(p_k) \ge C \log \frac1{g(|p_k|)}$;
\item
For all $n$, $m_{p\cup b} (n) \ge 1/ \tg(n)$.
\end{enumerate}
\end{lemma}

\begin{proof*}{\it Proof of Lemma \ref{mucantor} :}

This proof is patterned after that in \cite[pp. 51--52]{Lu-Se}, where 
instead of Conclusion 2, only the fact that $p$ was not Blaschke was 
to be obtained.

We define points $p_{n,j}$ associated to the arcs $I_{n,j}$:
\[
p_{n,j} := (1-2^{-n})\exp(2^{-n}{2\pi}(j+\frac12) i) 
\]
(their radial projection is at the center of the corresponding interval).

We will define the sequence through the choice of a certain
subfamily $\calf$ of the family of all the $I_{n,j}$. We shall define
simultaneously a sequence of probability measures $\mu_n$, supported
on the union of all the $I_{n,j}$ for a given index $n$, with 
the same uniform
density on each of the invervals where it is supported.

First we define a sequence of integers by
\begin{equation}
\label{defl}
l_n := \max \{ k \in \Zz_+ : k \le \log_2 \tilde g(n-j) + j, 
0\le j \le n \} .
\end{equation}

One verifies that $l_n \le l_{n+1} \le l_n +1$, and in fact the
sequence we have defined is the largest verifying that property
together with $l_n \le \log_2 \tg(n)$, for any nonnegative $n$.

The family $\calf_0:= \{ I_{n,j} : n \ge 0, j \in J_n\}$,
where $J_n$ is a subset of $\{0, \dots, 2^n -1\}$
which is defined recursively in the following way : if $l_{n+1} = l_n$,
we want 
\[
\bigcup \{ I_{n+1,j} : j \in J_{n+1} \} =
\bigcup \{ I_{n,j} : j \in J_{n} \} , 
\]
which is ensured by 
$J_{n+1} := \{ 2j, 2j+1 : j \in J_n \}$ ; if on the contrary
$l_{n+1} = l_n +1$, we only select the first half of each interval
at the $n$-th level, i.e. $J_{n+1} := \{ 2j : j \in J_n \}$ and
\begin{equation}
\label{jump}
\left| \bigcup \{ I_{n+1,j} : j \in J_{n+1} \} \right| =
\frac12 \left| \bigcup \{ I_{n,j} : j \in J_{n} \} \right| . 
\end{equation}
We will denote $p^0 := \{ p_{n,j} : I_{n,j} \in \calf_0 \}$.

By (\ref{jump}), we see that 
\begin{equation}
\label{nthlength}
\frac1{2\pi} \left| \bigcup \{ I_{n,j} : j \in J_{n} \} \right| 
=2^{-n} \#J_n = 2^{-l_n},
\end{equation}
and we thus define the probability measure $\mu_n$ as having
uniform density $2^{l_n}/2\pi$ on the set 
$\bigcup \{ I_{n,j} : j \in J_{n} \}$.

Observe that the way we have defined the families $J_n$ implies
that 
\[
\mu_m (I_{n,j}) = \mu_n (I_{n,j}), \forall j \in J_n, m \ge n ,
\]
and therefore if we denote by $\mu$ a weak limit point of the 
sequence $\{ \mu_n \}$, $\mu (I_{n,j}) = \mu_n (I_{n,j})= 2^{l_n - n}$. 

We let $h_0$ be the Poisson integral of the measure $\mu$. This is
a positive harmonic function with verifies 
\begin{equation}
\label{majorh}
h_0 (p_{n,j}) \ge c \frac{\mu_n (I_{n,j})}{|I_{n,j}|} = c 2^{l_n} ,
\forall I_{n,j} \in \calf_0.
\end{equation}

In order to ensure that our sequence satisfy Conclusion (1), we need to 
reduce the sequence $p^0$. First,
define a subset of $\calf_0$ by
$\calf := \{ I_{n,j} : j \in J_n , l_n \ge \log_2 \tg(n) -1 \}$,
then a subset $p$ of $p^0$ by
\begin{equation}
\label{equivp}
p := \{ p_{n,j} : I_{n,j} \in \calf \},
\mbox{ thus }2^{l_n} \le \tg (n) \le 2^{l_n+1} 
\mbox{ for } p_{n,j} \in p.
\end{equation}
Equations (\ref{equivp}) 
and (\ref{majorh}) imply that for an appropriate positive constant $C$,
$h:=Ch_0 (p_{n,j}) \ge \tg(n)$, for any $p_{n,j}$ in $p$, which
yields conclusion (1). 

Let  $b:= p^0 \setminus p$, where $p^0$ was defined after
equation (\ref{jump}).
Let us check conclusion (2). 

Since for any $(n,j)$, $j \in J_n$, there
exists  $j_k \in J_{k}$, $0\le k \le n-1$,
 such that 
\[
I_{n,j} \subset I_{n-1,j_{n-1}} \subset \dots \subset I_{0,j_0},
\]
we see that if $\alpha $ (the aperture of the Stolz angles) is
larger than an absolute constant, then $\phi_{p^0} (x) \ge n$
for any $x$ in $I_{n,j}$. (To deal with smaller values of $\al$,
we'd have to "thicken" the sequence by putting $N$ equidistant
points in each interval $I_{n,j} $, $N$ depending on $\al$.
This has no ill effects on the properties we're interested in.
Details are left to the reader). Then we have
\[
m_{p^0}(n) \ge \left| \bigcup \{ I_{n,j} : j \in J_{n} \} \right| =
2\pi 2^{-l_n}
\ge \frac1{\tg (n)},
\]
by the definition of $l_n$ and (\ref{nthlength}). 

\begin{claim}
The sequence $b=p^0 \setminus p$ is Blaschke.
\end{claim}

\begin{proof}
Let $A:= \{n: l_n < \log_2 \tg(n) -1 \}$; then 
\[
p^0 \setminus p = \{ p_{n,j} : j \in J_n , n \in A \} .
\]
We will show that if $n \neq m$ and $n$, $m \in A$, 
then $l_n \neq l_m$. Accepting this,
we then have, using (\ref{nthlength}),
\[
\sum_{p^0 \setminus p} (1-|p_k|) =
\sum_{n \in A} 2^{-n} \# J_n 
=
\frac1{2\pi} \sum_{n \in A} 2^{-l_n}
\le
\frac1{2\pi} \sum_k 2^{-k} < \infty, \mbox{ q. e. d.}
\]
To prove our first assertion :
if $n \in A$, then 
\[
l_n +1 < \log_2 \tg(n) \le \log_2 \tg(n+1) ,
\]
and (\ref{defl}) implies that, for any $n$ and $0 \le j \le n$,
\[
l_n + 1 \le \log_2 \tg(j) + (n-j) +1,
\]
therefore by the definition of $l_{n+1}$ as an upper bound,
$l_{n+1} \ge l_n +1$. This implies the desired property.
\end{proof}

If the sequence $p$ is empty, then by the above $l_{n+1} = l_ n +1$,  for any
$n$, so $l_ n=n +l_0$, and
 $$ 
\sum_n 1/\tilde{g}(n) \leq C \sum_n 2^{-n}<\infty ,
$$
in contradiction with our assumption.
\end{proof*}

\begin{proof*}{\it Proof of necessity in Theorem \ref{main} (\ref{sum}) :}

Given a function $g$ so that $\sum_0^\infty \frac{w_n}{\tilde g(n)}=\infty$,
consider the sequence $p$ given by Lemma \ref{mucantor}. Conclusion (2)
of the Lemma
then implies that $p \cup b \in \cals_w$, and Lemma \ref{union} then yields
that $p \in \cals_w$. Conclusion (1) and Lemma \ref{harm}
then show that $g$ is
not  an essential minorant for $\cals_w$.
\end{proof*}

\begin{proof*}{\it Proof of  necessity in Theorem \ref{main} (\ref{liminf}) :}

Suppose that $g$ is such that there exists a positive integer $C$ and a subset $E$
 of $\Zz_+$ such that $\sum_{n \in E} v_n < \infty$ and
\[
\sup_{n \notin E} \tilde g([\frac{n}C]) v_{n} =: A < \infty .
\]

Without loss of generality, we may assume
that $v_n \le 1$ for all $n$,
and apply Lemma \ref{mucantor} to the unique function
$g_1$ such that $g_1$ is constant on the intervals of 
the form $[1-2^{-n},1-2^{-n-1})$
and 
\[
\tg_1 (m) =\frac{A}{v_{C(m+1)-1}} = \max \{ \frac{A}{v_n} : [\frac{n}C] =m \} .
\]
  We get sequences $p$ and $b$ such that $b$ is Blaschke and
$m_{p\cup b} (n) \ge v_{C(n+1)}/A$, and a harmonic function $h$. 

Let $E_1 := 
\{ m = [\frac{n}C] \mbox{ for some }n \in E\}$.
For any $p_{m,j}$ in $p$, $m \notin E_1$,
$h(p_{m,j} ) \ge \tg_1 (m) \ge \tg(m) $. Define 
\[
a := \{ p_{m,j}  : m \in E_1 \}, \quad  q:= p \setminus a.
\]
Thus $h(p_{m,j} ) \ge \tg(m) $ whenever $p_{m,j} \in q$.
We still
need to prove that $a$ is Blaschke.
For any $m$ such that $p_{m,j} \in p$, by (\ref{nthlength})
and (\ref{equivp}),
\[
\sum_j (1-|p_{m,j} |) \asymp 2^{-l_m} \asymp \frac1{\tg_1(m)}
=\frac{v_{C(m+1)-1}}{A},
\]
so that
\[
\sum_{p_{m,j} \in a} (1-|p_{m,j} |) \preceq \sum_{n \in E} v_n < \infty,
\]
and $a$ is a Blaschke sequence, therefore $a \cup b$ is Blaschke too. 

Let $\tilde q$, $\tilde b$, $\tilde a$, $\tilde p$ be the sequences 
obtained respectively from $q$, $b$, $a$, $p$ by adjoining to each
point of the sequence a set of $M$ separated points located in a 
hyperbolic neighborhood of fixed size so that there exists $C'$ greater 
than $C$
such that
\[
m_{\tilde p \cup \tilde b} (C'n) \ge m_{p \cup b} (n) \ge v_{C(n+1)}/A .
\]

It is easy to see that $\tilde a \cup \tilde b$ is a Blaschke sequence
and, using Harnack's inequality, that there exists a positive constant $C_1$ such that 
$h(\tilde q_k) \ge C_1 \log (1/g(\tilde q_k))$,
for any $k$.
Now we want to show that $\tilde p \cup \tilde b = 
\tilde q \cup (\tilde a \cup \tilde b)
\in \call_v$. Take any $m$, and an integer $n$ such that 
$C(n+1) \le m \le C(n+2)$. For $m$ large enough, $m \le C(n+2) < C'n$, so
\[
m_{\tilde p \cup \tilde b} (m) \ge m_{\tilde p \cup \tilde b} (C'n)
\ge v_{C(n+1)}/A \ge v_m/A,
\]
therefore $q \cup (\tilde a \cup \tilde b) \in \call_v$, and since $\tilde a 
\cup \tilde b$ is Blaschke, $\tilde q \in \call'_v$. Lemma \ref{harm} then 
shows that $g$ is not an essential minorant for $\call_v$. 
\end{proof*}

\section{Comparison with summatory conditions}
\label{summatory}

The classes of sequences we have defined may seem somewhat arbitrary,
and perhaps one would like to express results about the decrease of 
holomorphic functions in terms of classes given by more usual
summatory conditions such as
\[
\sum_k f (1-|a_k|) = \infty,
\]
where $f$ is a positive increasing function on the interval $(0,1]$ such
that $\lim_{x\to0} f(x)/x$ $ =0$ (so that the condition is stronger than being
non-Blaschke), and $\int_0^1 x^{-2} f(x) dx = \infty$ (so that the condition
can be satisfied by some separated sequences).

There is no loss of generality in supposing that $f$ is constant on
intervals of the form $(2^{-n-1}, 2^{-n}]$, and if we set 
\[ 
w(n):=w_n := \frac{f(2^{-n})}{2^{-n}},
\]
then the above conditions become exactly those we have imposed on
$\swn$ in Definition \ref{classes}. Therefore, rather than
reasoning in terms of a function $f$, for a sequence $\{w_n\}$ as in 
Definition \ref{classes}, we define 
\begin{equation}
\label{sumw}
\calp_w := \left\{ a \mbox{ separated sequence }
: \sum_k (1-|a_k|) w([\log_2 \frac1{1-|a_k|}]) = \infty \right\}.
\end{equation}
Clearly this class is stable under removal of Blaschke sequences.

\begin{lemma}
\label{inclusion}
$\calp_w \subset \cals_w$, but for any $v$, and $w$
such that $\lim_{n\to\infty} w_n =0$, $\cals_v \not \subset \calp_w$.
\end{lemma}

On the other hand, it is also possible to determine the essential 
minorants for the class $\calp_w$.

\begin{theorem}
\label{thmsum}
$g$ is an essential minorant for the class $\calp_w$ if and only 
if 
$$
\sum_0^\infty \frac{w_n}{\tilde g(n)} < \infty .
$$
\end{theorem}

So the classes $\calp_w$ and $\cals_w$ have the same essential minorants,
while the former is quite a bit narrower than the latter. Since we gain 
information about the decrease of bounded holomorphic functions every 
time we can exhibit a sequence of points and a function $g$ which is an 
$H^\infty$-essential minorant for that sequence, Theorem
\ref{main}, by providing 
more sequences admitting a given essential minorant, seems a more 
interesting generalization of \cite{Lu-Se} than Theorem \ref{thmsum}.
Of course, it remains an open problem (and perhaps one which cannot admit 
any manageable answer) to determine the essential minorants for 
$H^\infty$ over a given sequence in the disk (rather than a class of 
sequences).

\begin{proof*}{\it Proof of Theorem \ref{thmsum} :}

The proof follows exactly the arguments in \cite{Lu-Se},
so we only sketch it. First, the following "weighted" version of 
\cite[Theorem 2]{Lu-Se} can be proved in exactly the same way.
For a real-valued function $u$ on the disk,
set $E_g (u) := \{ z \in \Dd : u(z) > \log \frac1{g(|z|)} \}$. We denote by
$\lambda_2$ the $2$-dimensional Lebesgue measure.

\begin{proposition}
\label{area}
For any $g$ as in Definition \ref{fcng}, $f$ as above, if
\[
\int_0^1 \frac{f(1-r) \, dr}{(1-r)^2 \log \frac1{g(r)}} < \infty,
\]
then for any (super)harmonic function $u$ on $\Dd$,
\[
\int_{E_g (u)} \frac{f(1-|z|) d\lambda_2 (z)}{(1-|z|)^2} < \infty.
\]
\end{proposition}

Proposition \ref{area} implies, as in \cite[pp. 50--51]{Lu-Se},
that the given condition is sufficient for $g$ to be an
essential minorant for $\calp_w$ (we actually have an if and only if 
statement, but don't need it right now).

Conversely, if we take $g$ such that 
$\sum_0^\infty \frac{w_n}{\tilde g(n)}= \infty$, Lemma \ref{mucantor}
yields a Blaschke sequence $b$, a separated sequence $p$ and a harmonic
function $h$ so that $h(p_k) \ge \log \frac1{g(1-|p_k|)}$. We also have,
for the sequence $p \cup b$,
\[
\sum_{z \in p \cup b, 1-|z|=2^{-n} } 1-|z| \ge \frac{C}{\tg(n)},
\]
so that 
\begin{multline*}
\sum_n w_n \left( \sum_{k: 1-|p_k| = 2^{-n} } 1-|p_k| \right) \\
\ge 
\sum_n w_n \left(  \frac{C}{\tg(n)} - \sum_{1-|b_k|=2^{-n}} 1-|b_k| \right)
\ge C \sum_0^\infty \frac{w_n}{\tilde g(n)} - \sum_k 1-|b_k|= \infty,
\end{multline*}
and therefore $p \in \calp_w$.
\end{proof*}

\begin{proof*}{\it Proof of Lemma \ref{inclusion} :}

Let $a \notin \cals_w$, $a$ separated. We want to prove that
$a \notin \calp_w$, that is to bound the sum in (\ref{sumw}).

Associate to any point $z$ in $\di$ the arc $I_z := \{\eit : z \in \Stolat \}$,
and define 
$$
W_a (\eit ) := \sum_k \chi_k w([\log_2 \frac1{1-|a_k|}]) ,
$$
where $\chi_k$ stands for the characteristic function of the arc
$I_{a_k}$. Since $|I_{a_k}| \asymp 1-| a_{k}|$,
the sum which we are trying to control is
bounded by a constant multiple of $\int_{\partial \di} W_a$.

Consider any point $\eit$ where $\phi_a (\eit)=n$. Then
$W_a (\eit ) = \sum_{j=1}^n w_{i_j}$, where 
$\Stolat \cap a = \{ a_{k_1}, \dots,  a_{k_n} \}$, and
$i_j = \left[ \log_2 \frac1{1-|a_{k_j}|} \right]$.
As in the proof
of Lemma \ref{phip}, we know that in 
$\Stolat \cap \{ 2^{-m} \ge 1-|z| > 2^{-m-1} \}$, there can be
no more than $M$ points of $a$; furthermore the monotonicity of
$w$ shows that the sum defining $W_a$ is largest when we take
the points $a_{k_j}$ with the smallest possible moduli. We thus get
$$
W_a (\eit )  \le  M \sum_{i=1}^{n/M} w_i \le  M \sum_{i=1}^n w_i .
$$

Recalling the definition of $m_a(n)$, we conclude with an integration by
parts:
$$
\int_{\partial \di} W_a (\theta) d \theta  \le
 M \sum_{n \ge 1} \left( \sum_{i=1}^n w_i  \right) (m_n - m_{n+1})
\le 
\sum_n m_a(n) w_n < \infty .
$$

To see that $\calp_w$ contains no $\cals_v$ when $w_n$ tends to $0$,
consider a sequence of the form 
\[
a_{k,j} := (1-2^{-n_k}) \exp(2 \pi i j 2^{-n_k}), \quad
0 \le j < 2^{n_k}.
\]
Picking $\{ n_k \}$ an increasing sequence of integers so
that $\sum_k w(n_k) < \infty$, we see that $a \notin \calp_w$; on the
other hand, it is easy to see that $\phi_a (\eit) = \infty$ almost 
everywhere on $\partial \Dd$ (for an $\al$ large enough, and therefore 
for any $\al$). Thus $a \in \cals_v$, for any $v$ such that $\sum_n v_n = 
\infty$.
\end{proof*}

We'd like to thank our colleague Artur Nicolau for useful suggestions. The
second-named author also would like 
to thank Boris Korenblum for his hospitality and
stimulating conversations.

\vskip1cm

Jordi Pau, Departament de Matem\`atiques, Universitat Aut\`onoma
de Bar\-ce\-lo\-na, 08193 Bellaterra, Spain.

Pascal J. Thomas, Laboratoire de Math\'ematiques Emile Picard, UMR CNRS 5580,
Universit\'e Paul Sabatier, 118 route de Narbonne, 31062 TOULOUSE CEDEX, France.

                    \end{document}